\documentclass[twoside, letterpaper, 11pt]{amsart}
\usepackage{amsfonts,amsmath,amssymb,amsthm,enumerate}
\usepackage{eucal}
\usepackage{graphicx}
\usepackage[usenames,dvipsnames,svgnames,table]{xcolor}
\graphicspath{{./figures/}}
\usepackage{hyperref}
\hypersetup{colorlinks,
filecolor=black,
linkcolor=MidnightBlue,
citecolor=NavyBlue,
urlcolor=RoyalBlue,
bookmarksopen=true}
\usepackage{lineno}

\usepackage{float}
\allowdisplaybreaks

\usepackage[letterpaper, total={6.3in, 8in}, left=1.1in, top=1.5in]{geometry}

\usepackage[alphabetic]{amsrefs}
\AtBeginDocument{%
	\def\MR#1{}
}

\theoremstyle{plain}

\theoremstyle{remark}

\newtheorem*{remark}{Remark}

\newcommand{\cv}{\nabla}

\newcommand{\mb}{\mathbb}
\newcommand{\mc}{\mathcal}

\newcommand{\mr}{\mathrm}
\newcommand{\ten}{\otimes}

\newcommand{\ve}{\varepsilon}
\newcommand{\vp}{\varphi}

\newcommand{\mat}{\begin{pmatrix}}
\newcommand{\rix}{\end{pmatrix}}

\DeclareMathOperator{\Rc}{Rc}
\newcommand{\icpn}{\int_{\mb{CP}^2}}

\begin{document}

\title[Asymptotics of unstable perturbations of the Fubini--Study metric in Ricci flow]{Asymptotic behavior of unstable perturbations of the Fubini--Study metric in Ricci flow}

\author{David Garfinkle}
\address[David Garfinkle]{Oakland University}
\email{garfinkl@oakland.edu}

\author{James Isenberg}
\address[James Isenberg]{University of Oregon}
\email{isenberg@uoregon.edu}
\urladdr{http://www.uoregon.edu/$\sim$isenberg/}

\author{Dan Knopf}
\address[Dan Knopf]{University of Texas at Austin}
\email{danknopf@math.utexas.edu}
\urladdr{http://www.ma.utexas.edu/users/danknopf}

\author{Haotian Wu}
\address[Haotian Wu]{The University of Sydney}
\email{haotian.wu@sydney.edu.au}

\thanks{
DG thanks the NSF for support in PHY-1806219 and PHY-2102914. 
JI thanks the NSF for support in PHY-1707427.
DK thanks the Simons Foundation for support in Award 635293.
HW thanks the Australian Research Council for support in DE180101348.
}

\begin{abstract}
Kr\"oncke has shown that the Fubini--Study metric is an unstable generalized stationary solution of Ricci flow~\cite{Kro20}. In this paper, we carry out
numerical simulations which indicate that Ricci flow solutions originating at unstable perturbations of the Fubini--Study metric develop local singularities modeled by
the blowdown soliton discovered in~\cite{FIK03}.
\end{abstract}

\maketitle

\tableofcontents

\section{Introduction}

In studying dynamical systems, it is always important to identify and classify fixed points. For Ricci flow, regarded as a dynamical system on the
infinite-dimensional space of Riemannian metrics, one must consider generalized fixed points. These are the \emph{Ricci solitons}, specified by the data $\big(\mc M^n,g,\lambda, X\big)$ for which
\[
-2\Rc[g]  = 2\lambda g + \mc L_X g,
\]
where $\Rc[g]$ is the Ricci curvature of the metric $g$, $\lambda\in\mathbb{R}$, and $\mc L_X g$ is the Lie derivative of $g$ with respect to the vector field $X$. 

A prominent example of a Ricci soliton is the Fubini--Study metric on complex projective space, an Einstein metric corresponding to $\big(\mb{CP}^N,G_{\mr{FS}},2(N+1),0\big)$.\footnote{
$N$ is the complex dimension, so that the real dimension is $n=2N$.}
In \S10 of~\cite{Ham95}, Hamilton conjectures that it is a stable fixed point. We refer the reader to Section \ref{initial-data} for an explicit description of the Fubini--Study metric. In Example~2.3 of~\cite{CHI04}, Cao, Hamilton, and Ilmanen note that the Fubini--Study
metric is neutrally linearly stable with respect to the second variation of Perelman's shrinker entropy, crediting this observation to unpublished work of Goldschmidt. However, in his dissertation,  Kr\"oncke computes the third variation of Perelman's entropy to prove the surprising result that Fubini--Study is in fact unstable. This calculation has appeared as~\cite{Kro20} and has been independently verified in~\cite{KS19}. The instability is conformal and hence not K\"ahler.

The instability of the Fubini--Study raises an interesting question. What happens to Ricci flow solutions that start at arbitrarily small but unstable perturbations of~$G_{\mr{FS}}$? According to the conjectural hierarchy of $4$-dimensional solutions given in~\cite{CHI04}, one expects the Ricci flow of unstable perturbations of $\big(\mb{CP}^2,G_{\mr{FS}}\big)$ to develop local singularities modeled by the blowdown gradient shrinking soliton $\mc L^2_{-1}$ discovered in~\cite{FIK03}. These solutions would start arbitrarily close to (but \emph{not} in) the subset of K\"ahler metrics and, in finite time, develop local singularities that are asymptotically K\"ahler (but with the complex structure reversed --- see below). A closely related type of singularity formation, an asymptotic approach to the space of K\"ahler metrics in the blowup limit, has been partially investigated in~\cite{IKS19}. It should be noted that, while the blowdown  soliton has been shown to occur as a model of singularity formation on compact manifolds by M\'aximo~\cite{Max14}, the solutions considered in that work are K\"ahler, whereas the unstable perturbations studied here all lie outside the space of K\"ahler initial data.

In this paper, we describe numerical simulations that provide evidence in favor of the conjecture that there exists an unstable Ricci flow ``orbit'' starting arbitrarily close to $\big(\mb{CP}^2, G_{\mr{FS}}\big)$ and ending at the blowdown soliton $\mc L^2_{-1}$ in real dimension $n=4$. This is not a true orbit, because convergence to the blowdown soliton happens only locally and after parabolic dilation of the developing singularity. Nonetheless, it provides evidence in favor
of the hierarchy outlined in~\cite{CHI04}. Furthermore, while $\mc L^2_{-1}$ is known to be unique among K\"ahler solitons, our results suggest that there is no
non-K\"ahler gradient shrinking soliton with the topology of $\mb C^2$ blown up at a point, at least none with a density $\Theta$ above that of $\mb{CP}^2$; \emph{cf}.~\cite{CHI04}.
\medskip

Because unperturbed $\big(\mb{CP}^2,G_{\mr{FS}}\big)$ is a positively-curved Einstein manifold, under Ricci flow it vanishes in a point. In particular, the area of the distinguished $\mb{CP}^1$ goes to zero at the same time that the diameter of the manifold goes to zero. In each of our numerical simulations, we begin with initial data that are small non-K\"ahler conformal perturbations $\big(\mb{CP}^2,\tilde G\big)$. The conformal factor lies in the nullspace of the second variation of Perelman's
shrinker entropy $\nu$ and is constructed in Section~\ref{unstable}. Our numerics indicate that evolution of this initial data by Ricci--DeTurck flow crushes the
$\mb{CP}^1$  fiber before the remainder of the manifold vanishes, that is, before the diameter goes to zero. Performing a parabolic rescaling of the solution fixes the
area of the $\mb{CP}^1$ fiber. Our numerics indicate that the developing singularity is Type-I; that is, it develops at the natural parabolic rate. Because of this, one expects by~\cite{EMT11} to see local convergence of the parabolic blowups to a nonflat gradient shrinking soliton. How do we recognize the soliton numerically?
We do so in two ways.

First, using a measure of local closeness to K\"ahler derived in Section~\ref{FIK soliton} below, we find that the blowups are becoming asymptotically K\"ahler near the
distinguished fiber, but with the opposite orientation of the complex structure on unperturbed $\big(\mb{CP}^2,G_{\mr{FS}}\big)$. This is significant because we know
by~\cite{FIK03} that the blowdown soliton is the unique gradient shrinking K\"ahler soliton with the topology of $\mb C^2$ blown up at the origin, which corresponds
in the original solution to the distinguished $\mb{CP}^1$.

Second, we use the fact that the blowdown soliton is asymptotically conical. By~\cite{KW15}, there is at most one shrinking gradient soliton
asymptotic to any cone. So, as further verification, we determine that the parabolic blowup is converging numerically to the correct cone  near the $\mb{CP}^1$,
in a neighborhood suitable to the rescaling. Thus, while numerical simulations cannot provide proof, we regard these two observations as compelling evidence in favor of the formation of the~\cite{FIK03} soliton in solutions emerging from unstable perturbations of the Fubini--Study metric.

\section{An unstable conformal perturbation}		\label{unstable}

We seek a conformal factor $\psi$ so that $(1+\delta\psi)G_{\mr{FS}}$ is an unstable perturbation of the Fubini--Study metric for small $\delta>0$.

As shown in~\cite{KS19}, the function $\psi$ must satisfy three conditions:
\[
\Big(\Delta+\frac{1}{\tau}\Big)\psi=0,\qquad\icpn\psi\,\mr dV=0,\qquad \icpn\psi^3\,\mr dV>0.
\]
These are, respectively: the condition that $\psi$ lie in the nullspace of the second variation of
Perelman's entropy $\nu$, a normalization necessary for that functional, and a sufficient condition for the perturbation generated by $\psi$ to be unstable.
Here, $1/(2\tau)$ is the Einstein constant of the metric, which in our case is $2(N+1)=6$.

Because $\Rc(G_{\mr{FS}})=6\,G_{\mr{FS}}$ and the Laplacian of a rotationally-symmetric function $\psi$ is 
\[
\Delta \psi = \psi_{\theta\theta}+\big\{3\cot(\theta)-\tan(\theta)\big\}\psi_\theta,
\]
to specify rotationally symmetric unstable perturbations of the Fubini--Study metric, we seek solutions $\psi$ of the second-order linear \textsc{ode}
\begin{align*}
0	&=(\Delta+12)\psi(\theta)\\
	&=\psi''(\theta)+\big\{3\cot(\theta)-\tan(\theta)\big\}\psi'(\theta)+12\psi(\theta).
\end{align*}
Only one fundamental solution of this \textsc{ode} is smooth for $\theta\in[0,\pi/2]$. Up to an arbitrary multiplicative constant, it is
\[
\psi=\csc^2(\theta)\big\{2-8\cos(2\theta)+6\cos(4\theta)\big\} =1+3\cos(2\theta).
\]
Using $\det(G_{\mr{FS}})=\sin^6(\theta)\cos^2(\theta)$, one readily verifies that this solution satisfies
\begin{align*}
\icpn\psi\,\mr dV = \int_0^{\pi/2}\big\{1+3\cos(2\theta)\big\}\sin^3(\theta)\cos(\theta)\,\mr d\theta=0,\\
\icpn\psi^3\,\mr dV = \int_0^{\pi/2}\big\{1+3\cos(2\theta)\big\}^3\sin^3(\theta)\cos(\theta)\,\mr d\theta=\frac25.
\end{align*}

Hence, given $\delta>0$, we choose as initial data
\[
\tilde G = \frac{1+\delta\psi}{1+\delta}G_{\mr{FS}}=\Big(1+3\frac{\delta}{1+\delta}\cos(2\theta)\Big)G_{\mr{FS}}=h\,G_{\mr{FS}},
\]
where $h$ is defined in~\eqref{define h} below, with $\ve=\frac{\delta}{1+\delta}$.

\section{Cohomogeneity-one metrics}	\label{setup}

On $[0,\pi/2]\times\mr{SU}(2)$, we consider cohomogeneity-one metrics of the form
\begin{equation}	\label{metric}
G=\mr \rho^2\,d\theta^2+\sum_{i=1}^3 f_i^2\,\omega^i\ten\omega^i,
\end{equation}
where $\omega^i\ten\omega^i$ is a Milnor coframe for $\mr{SU}(2)$.
The cohomogeneity-one condition means that the fiber over a generic $\theta$ is diffeomorphic to $\mb S^3$, except for $\theta\in\{0,\pi/2\}$.
We are interested in initial data that are unstable conformal perturbations of the Fubini--Study metric on $\mb{CP}^2$. Accordingly,
we impose an \emph{Ansatz} of $\mr U(2)$ symmetry in the form
\begin{equation}	\label{symmetry}
f_1=f\qquad\mbox{ and }\qquad f_2=f_3=g.
\end{equation}

From the viewpoint of analysis, it is convenient to fix a gauge and let $s(\theta,t)$ represent arclength with respect to $G$. Then $s$ and $\theta$ satisfy the infinitesimal spatial relation
\[
\mr ds = \rho\,\mr d\theta.
\]
In this gauge, the metric takes the form
\[G= ds^2 +\sum_{i=1}^3 f_i^2\,\omega^i\ten\omega^i.\]
Then the sectional curvatures of $G$ are convex linear combinations of
\begin{subequations}		\label{sectionals}
	\begin{align}
		\kappa_{12}=\kappa_{31}&=\frac{f^2}{g^4}-\frac{f_s g_s}{fg},\\
		\kappa_{23}&=\frac{4g^2-3f^2}{g^4}-\frac{g_s^2}{g^2},\\
		\kappa_{01}&=-\frac{f_{ss}}{f},\\
		\kappa_{02}=\kappa_{03}&=-\frac{g_{ss}}{g},
	\end{align}
\end{subequations}
where $e_0=\frac{\mr d}{\mr ds}$ and $(e_1,e_2,e_3)$ is a Milnor frame for $\mr{SU}(2)$.

The variables $(s,t)$ do not commute, but this approach yields a strictly parabolic system for Ricci flow on these geometries,
\begin{subequations}		\label{parabolic}
	\begin{align}
		f_t		&=f_{ss}+2\frac{g_s}{g}f_s-2\frac{f^3}{g^4},\\
		g_t	&=g_{ss}+\left(\frac{f_s}{f}+\frac{g_s}{g}\right)g_s+2\frac{f^2-2g^2}{g^3},
	\end{align}
\end{subequations}
where all time derivatives on the \textsc{lhs} are taken at fixed $\theta$. This system has been studied in~\cite{IKS19}. One can recover the behavior of $\rho$ under Ricci flow using its evolution equation
\[\rho_t = \left(\frac{f_{ss}}{f}+2\frac{g_{ss}}{g}\right)\rho.\]

From the viewpoint of numerical simulation, on the other hand, it is preferable to use fixed commuting variables.
Using the fact that any smooth $\zeta(s)$ must satisfy
\[
\zeta=\frac{\zeta_\theta}{\rho}\qquad\mbox{ and }\qquad\zeta_{ss}=\frac{\zeta_{\theta\theta}}{\rho^2}-\frac{\rho_\theta\zeta_\theta}{\rho^3},
\]
we rewrite system~\eqref{parabolic}  in terms of the fixed commuting variables $(\theta,t)$ as
\begin{subequations}		\label{NotParabolic}
\begin{align}
\rho_t&= \frac{f_{\theta\theta}}{\rho f}-\frac{\rho_\theta f_\theta}{\rho^2 f}+2\frac{g_{\theta\theta}}{\rho g}-2\frac{\rho_\theta g_\theta}{\rho^2g}, \label{dtrho}\\
f_t	&= \frac{f_{\theta\theta}}{\rho^2}-\frac{\rho_\theta f_\theta}{\rho^3}+2\frac{f_\theta g_\theta}{\rho^2g}-2\frac{f^3}{g^4},\\
g_t	&=\frac{g_{\theta\theta}}{\rho^2}-\frac{\rho_\theta g_\theta}{\rho^3}+
	\left(\frac{f_\theta}{f}+\frac{g_\theta}{g}\right)\frac{g_\theta}{\rho^2}+2\frac{f^2-2g^2}{g^3}.
\end{align}
\end{subequations}
The evolution equation \eqref{dtrho} for $\rho$ prevents this system from being parabolic. To remedy this, we implement DeTurck's trick below.

\section{Initial data} \label{initial-data}
Our reference (unperturbed) Fubini--Study metric $G_{\mr{FS}}$ is determined by three functions
$\rho_{\mr{FS}}, f_{\mr{FS}}, g_{\mr{FS}}\colon [0,\pi/2]\rightarrow\mb R$, given by
\[
\rho_{\mr{FS}}(\theta)=1,\qquad
f_{\mr{FS}}(\theta)=\frac12\sin(2\theta)=\sin(\theta)\cos(\theta),	\quad\mbox{ and }\quad g_{\mr{FS}}(\theta)=\sin(\theta).
\]
That is, using \eqref{metric}, we have
\[ G_{\mr{FS}} = d\theta^2 + \sin^2(\theta)\cos^2(\theta)\, \omega^1\ten\omega^1+\sin^2(\theta)\,\omega^2\ten\omega^2+\sin^2(\theta)\,\omega^3\ten\omega^3. \]

It is readily verified that these choices have the correct boundary behavior to induce a smooth metric on $\mb{CP}^2$, with the geometry of an
asymptotically round $4$-ball near $\theta=0$ and with a distinguished $\mb{CP}^1$ fiber at $\theta=\pi/2$.  Furthermore, it follows easily from~\eqref{sectionals}
that $\Rc[G_{\mr{FS}}]=6\,G_{\mr{FS}}$.

For use in specifying a Ricci-DeTurck flow below, which shadows the Ricci flow and is parabolic, unlike the Ricci flow itself, we note that the Fubini--Study connection is determined by
\begin{equation}
(\Gamma_{\mr{FS}})_{00}^0 = 0,\quad
(\Gamma_{\mr{FS}})_{11}^0 = -\frac14\sin(4\theta),\quad
(\Gamma_{\mr{FS}})_{22}^0=(\Gamma_{\mr{FS}})_{33}^0=-\frac12\sin(2\theta).
\end{equation}
\medskip

Given any $\ve>0$, we recall that the conformal factor determining an unstable perturbation of the Fubini--Study metric is given by
\begin{equation}		\label{define h}
h(\theta):=\big(1+3\ve\cos(2\theta)\big),
\end{equation}
as shown in Section~\ref{unstable}. 
As initial data for our perturbed Ricci flow, we take $G(0)=\tilde G=h\,G_{\mr{FS}}$, for which $(\rho,f,g)=(h\rho_{\mr{FS}},hf_{\mr{FS}},hg_{\mr{FS}})$, namely
\begin{align*}
\rho(\theta)	&=\Big(1+3\ve\cos(2\theta)\Big),\\
f(\theta)		&=\frac14\Big(2\sin(2\theta)+3\ve\sin(4\theta)\Big),\\
g(\theta)		&=\sin(\theta)\Big(1+3\ve\cos(2\theta)\Big)
\end{align*}
at time $t=0$. The Levi--Civita connection of $G(0)$ is determined by
\begin{align*}
\Gamma_{00}^0	&=-\frac{6\ve\sin(2\theta)}{1+3\ve\cos(2\theta)},\\
\Gamma_{11}^0	&=-\frac{\sin(2\theta)\Big(\cos(2\theta)+3\ve\cos(4\theta)\Big)}{2\Big(1+3\ve\cos(2\theta)\Big)},\\
\Gamma_{22}^0=\Gamma_{33}^0	&=-\frac{\sin(\theta)\cos(\theta)\Big(1-6\ve+9\ve\cos(2\theta)\Big)}{1+3\ve\cos(2\theta)},
\end{align*}
with the unlisted Levi-Civita coefficients vanishing.

\section{Ricci--DeTurck flow}
The formulas below that suppress time apply at all $t\geq0$ such that a solution exists.\medskip

We consider the \textsc{pde} system
\begin{equation}
\partial_t G = -2\Rc[G]+\mc L_VG
\end{equation}
with initial data $G(0)=h\,G_{\mr{FS}}$, where our DeTurck vector field $V$ is defined by
\begin{subequations}
\begin{equation}
V^\beta = G^{\alpha\alpha}\Big(\Gamma_{\alpha\alpha}^\beta-(\Gamma_{\mr{FS}})_{\alpha\alpha}^\beta\Big),
\end{equation}
in which the implicit summation ranges over $\alpha\in\{0,1,2,3\}$.
For the metrics under consideration, this reduces to $V=V^0\,\frac{\mr d}{\mr d\theta}$, where
\begin{equation}
V^0 = G^{\alpha\alpha}\Big(\Gamma_{\alpha\alpha}^0-(\Gamma_{\mr{FS}})_{\alpha\alpha}^0\Big).
\end{equation}
\end{subequations}

\begin{remark}
We note that at $t=0$, one has
\[
V^0(\theta,0) = \frac{12\ve\sin(2\theta)}{\Big(1+3\ve\cos(2\theta)\Big)^3},
\]
which is smooth for all sufficiently small values of $\ve\geq0$.
\end{remark}

In working with the DeTurck vector field $V=V^0\,\frac{\mr d}{\mr d\theta}$ and the 1-form $V^\flat=V_0\,\mr d\theta$, where $V_0 = \rho^2 V^0$,
it is convenient to relabel $V^0$ as the function $v(\theta,t)$ defined by
\begin{equation}		\label{DeTurck}
v= \frac{\rho_\theta}{\rho^3} + \frac{\tfrac14\rho^2\sin(4\theta)-ff_\theta}{\rho^2 f^2}
	+\frac{\rho^2\sin(2\theta)-2gg_\theta}{\rho^2 g^2}.
\end{equation}
One finds that $(\mc L_V G)_{\alpha\alpha} = 2\cv_\alpha V_\alpha$, where the only nonzero components are
\begin{align*}
\cv_0 V_0	&= \rho^2v_\theta+v\,\rho\rho_\theta = \rho(\rho v)_\theta,\\
\cv_1 V_1	&= v\,ff_\theta,\\
\cv_2 V_2	&= \cv_3 V_3 = v\,gg_\theta.
\end{align*}
Thus from system~\eqref{NotParabolic} above, we obtain a manifestly parabolic system corresponding to the Ricci-DeTurck flow,
\begin{subequations}		\label{DansForm}
\begin{align}
\rho_t&= \frac{\rho_{\theta\theta}}{\rho^2}-3\frac{\rho_\theta^2}{\rho^3}+\Big(\frac{f_\theta}{\rho^2f}+2\frac{g_\theta}{\rho^2g}+v\Big)\rho_\theta\\
	&\quad+\left( \Big[\frac{\sin(4\theta)}{4f^2}\Big]_\theta+\Big[\frac{\sin(2\theta)}{g^2}\Big]_\theta\right)\rho
	+\frac{f_\theta^2}{\rho f^2}+2\frac{g_\theta^2}{\rho g^2},	\notag \\
f_t	&= \frac{f_{\theta\theta}}{\rho^2}+\Big(v-\frac{\rho_\theta}{\rho^3}\Big)f_\theta+2\frac{f_\theta g_\theta}{\rho^2g}-2\frac{f^3}{g^4},\\
g_t	&=\frac{g_{\theta\theta}}{\rho^2}+\Big(v-\frac{\rho_\theta}{\rho^3}\Big)g_\theta+
	\left(\frac{f_\theta}{f}+\frac{g_\theta}{g}\right)\frac{g_\theta}{\rho^2}+2\frac{f^2-2g^2}{g^3}.
\end{align}
\end{subequations}

Below, we convert this into an equivalent form for which suitable boundary behaviors as $\theta\searrow0$ and $\theta\nearrow\pi/2$ can be robustly enforced
in numerical simulations.

\section{A system better suited to simulation}
Here, we convert system~\eqref{DansForm} into an equivalent one appropriate for numerical simulation of the evolution. The issue we must address is this: 
for the evolving functions $(\rho,f,g)$ to induce smooth metrics on $\mb{CP}^2$ at each time requires that we enforce the identities
\[
f=0,\quad f_s=1,\quad\mbox{ and }\quad g=0,\quad g_s=1 \qquad\mbox{ at }\qquad \theta=0,
\]
and the identities
\[
f=0,\quad f_s=-1,\quad\mbox{ and }\quad g>0,\quad g_s=0 \qquad\mbox{ at }\qquad \theta=\frac{\pi}{2}.
\]
However, we may impose only one Dirichlet or Neumann condition at each boundary for the functions that we evolve numerically. We proceed to resolve this issue in two steps, as follows:
\medskip

We begin by defining
\[
\tilde\rho:=\log\rho,\qquad\tilde f:=\log f,\qquad\tilde g:=\log g,
\]
and
\[
Q:=\rho^2\sin(\theta)\cos(\theta)\left(\frac{\cos(2\theta)}{f^2}+\frac{2}{g^2}\right).
\]
Then using~\eqref{DansForm} and diligently doing some algebra, we obtain the system:
\begin{subequations}		\label{DavidsForm}
\begin{align}
\rho^2\,\tilde\rho_t	&=\tilde\rho_{\theta\theta}-Q\tilde\rho_\theta-\tilde\rho_\theta^2+\tilde f_\theta^2+2\tilde g_\theta^2+ Q_\theta,\\
\rho^2\,\tilde f_t		&=\tilde f_{\theta\theta}+Q\tilde f_\theta-2\frac{\rho^2}{g^4}f^2,\\
\rho^2\,\tilde g_t	&=\tilde g_{\theta\theta}+Q\tilde g_\theta+\frac{\rho^2}{g^4}\big(2f^2-4g^2\big).
\end{align}
\end{subequations}

Next we define $A,Y,Z$ via the relations
\[
\rho=:e^{A+Y+Z},\qquad f=:e^Y\cos(\theta)\,g,\qquad g=:e^A\sin(\theta),
\]
noting that
\[
Q=2e^{2Z}\left(\cot(2\theta)+e^{2Y}\cot(\theta) \right).
\]
Then we obtain this system:
\begin{subequations}
\begin{align}
\rho^2\,A_t	&= A_{\theta\theta}+QA_\theta +\csc^2(\theta)\left( 2e^{4Y+2Z}\cos^2(\theta)-4e^{2Y+2Z}-1\right) \label{dtA}\\
	&\quad+2e^{2Z}\cot(\theta)\left(\cot(2\theta)+e^{2Y}\cot(\theta)\right),\notag\\
\rho^2\,Y_t	&=Y_{\theta\theta}+QY_\theta-2e^{2Z}\tan(\theta)\cot(2\theta)-2e^{2Y+2Z}-\sec^2(\theta) \label{dtY}\\
	&\quad+4e^{2Y+2Z}\csc^2(\theta)\left(1-e^{2Y}\cos^2(\theta)\right),\notag\\
\rho^2\,Z_t	&=Z_{\theta\theta}-Q(2A_\theta+2Y_\theta+Z_\theta)+Q_\theta-(A_\theta+Y_\theta+Z_\theta)^2+3\big(A_\theta+\cot(\theta)\big)^2 \label{dtZ}\\
		&+\big(Y_\theta-\tan(\theta)\big)^2+2\big(A_\theta+\cot(\theta)\big)\big(Y_\theta-\tan(\theta)\big)+2e^{2Y+2Z} \notag \\
		&+2e^{2Y+2Z}\cot^2(\theta)(e^{2Y}-1)+\sec^2(\theta)+\csc^2(\theta)+2e^{2Z}\cot(2\theta)\big(\tan(\theta)-\cot(\theta)\big). \notag
\end{align}
\end{subequations}

Now consider what boundary conditions we need to impose smoothness.  We need ${A_\theta}=0 $ at $\theta = 0$ and $\theta =\pi/2$. We need $Y=0$ and ${Y_\theta}=0 $ at $\theta = 0$ and ${Y_\theta}=0$ at $\theta = \pi/2$.  We need $Z=0$ and ${Z_\theta}=0$ at both $\theta =0$ and $\theta = \pi/2$.  Thus, $A$ is the sort of variable we need, but $Y$ and $Z$ are not.  However, defining $B$ and $C$ by 
\[Y =: B {\sin ^2} (\theta),\qquad Z =: C {\sin^2} (2 \theta),\] 
we see that the conditions needed for smoothness can be imposed as long as ${B_\theta}={C_\theta}=0$ at $\theta=\pi/2$.  

After some tedious algebra we find that equation \eqref{dtY} yields

\begin{align} \label{dtB}
{\partial _t} B &= {\rho ^{-2}} \biggl [ {B_{\theta \theta}} + \left(Q + 4 \cot (\theta)\right) {B_\theta} + B \left ( {{2 + 2 {e^{2Z}} -4 {e^{2Y+2Z}}} \over {{\sin^2} (\theta)}} \right ) 
\\
&- 4 B (1 + {e^{2Z}} + {e^{2Y+2Z}}) 
+ {e^{2Z}} \left ( {{-2 - 2 {e^{2Y}} + 4 {e^{4Y}}} \over {{\sin^2}(\theta)}} \right ) 
\notag
\\
&+ 4 \left ( {{e^{2Z} - 1} \over {{\sin^2} (2 \theta)}} \right ) + 4 {e^{2Y+2Z}} \left ( {{1+2Y-{e^{2Y}}} \over {{\sin^4}(\theta)}} \right ) \biggr ].
\notag
\end{align}

Similarly some tedious algebra applied to equation \eqref{dtZ} yields the following:
\begin{align} \label{dtC}
{\partial _t} C &= {\rho ^{-2}} \biggl [ {C_{\theta \theta}} +\left(Q+8\cot (2\theta)\right){C_\theta} 
-8 C (2 + {e^{2Z}} + {e^{2Y+2Z}}) 
\\ &+ \left ( {{1-{e^{2Y+2Z}}} \over {{\sin^2}(\theta)}} \right ) \left ( {{2{A_\theta}} \over {\sin (2\theta)}} - 4 C \right )
\notag
\\
&+ \left ( {{1-{e^{2Z}}} \over {{\sin^2}(2\theta)}} \right ) \left ( -8C + 4 \cot (2\theta) {A_\theta} + 2 \cos (2\theta) \left( \tan(\theta) {B_\theta} +2B\right) - 6 \right )
\notag
\\
&+ {e^{2Y+2Z}} \left ( {{\cosh (2 Y) - 1} \over {{\sin^4} (\theta)}} \right ) + 8 \left(1+{\cos^2} (\theta)\right) \left ( {{1+2Z -{e^{2Z}}} \over {{\sin^4} (2 \theta)}} \right )
+ 2{{\left ( {{A_\theta}\over {\sin (2 \theta)}} \right ) }^2} 
\notag
\\
&- \left(\sin (2\theta) {C_\theta} + 4 \cos (2 \theta) C \right) \left ( \sin (2\theta) {C_\theta} + 4 \cos (2 \theta) C  + \tan (\theta) {B_\theta} + 2 B + {{2{A_\theta}}\over {\sin (2 \theta)}} \right ) \biggr ].
\notag
\end{align}

In summary, we evolve the variables $A$, $B$ and $C$, and their evolution equations are equations \eqref{dtA}, \eqref{dtB}, and \eqref{dtC}, respectively.  The Neumann boundary conditions are that ${A_\theta}, \; {B_\theta}$ and $C_\theta$ vanish at $\theta=0$ and $\theta =\pi/2$.

\section{The blowdown soliton}		\label{FIK soliton}
We first derive a necessary and sufficient condition for a metric of the form~\eqref{metric} to be K\"ahler, under the symmetry assumptions $f_1=f$ and $f_2=f_3=g$.
Calabi observed \cite{Cal82} that any $\mr U(2)$-invariant Kähler metric on $\mb C^2\backslash(0,0)$  may be written in complex form as
\begin{equation} 	\label{eqn:U2 Complex}
G_{\mb C} =\Big\{e^{-r}\vp\,\delta_{\alpha\beta} +e^{-2r}(\vp_r-\vp)\,\bar z_\alpha  z_\beta\Big\}\, \mr dz^\alpha\ten\mr d\bar z^\beta,
\end{equation}
with respect to his coordinate $r:=\log(|z_1|^2+|z_2|^2)$. 
Written in real coordinates, the same metric becomes
\begin{equation} 	\label{eqn:U2 Real}
  G_{\mb R} =\vp_r\Big(\frac14\mr dr\ten\mr dr+\omega^1\ten\omega^1\Big)
  +\vp\Big(\omega^2\ten\omega^2+\omega^3\ten\omega^3\Big),
\end{equation}
which shows that $\vp>0$ and $\vp_r>0$ are necessary conditions for this to be a K\"ahler metric.

A comparison of equations~\eqref{metric} and~\eqref{eqn:U2 Real} shows that a coordinate transformation is needed to write a K\"ahler metric with
respect to arclength $s$, where $\mr ds = \rho\,\mr d\theta$. To see this, we observe that if $s$ and $r$ are related by the \textsc{ode}
\[
\frac{\mr dr}{\mr ds}=\frac{2}{f},
\] 
then equation~\eqref{metric} takes the form
\begin{equation}		\label{eqn:metric-r}
  G = f^2\Big(\frac14 \mr d r\ten\mr d r+\omega^1\ten\omega^1\Big) + g^2\Big(\omega^2\ten\omega^2+ \omega^3\ten\omega^3\Big),
\end{equation}
which matches~\eqref{eqn:U2 Real} if and only if $f^2=\vp_r$ and $g^2=\vp$ are related by
\begin{equation} 	\label{eqn:KahlerCondition}
f=gg_s.
\end{equation}
It follows that a cohomogeneity-one metric with the $\mr U(2)$ symmetries we have imposed is K\"ahler if and only if condition~\eqref{eqn:KahlerCondition} holds.
So, we may use the ratio $f/(gg_s)$ to provide a local measurement of the closeness of a metric to the K\"ahler subspace.
\medskip

Now, it follows from Lemma~6.1 and equation~(27) of ~\cite{FIK03} with $\lambda=-1$, $\mu=\sqrt2$, and $\nu=0$, that the metric on the blowdown
soliton $\mc L^2_{-1}$  is determined by a function $\phi(r)$ that solves the separable first-order \textsc{ode}
\begin{equation}	\label{FIK ode}
 \phi_r=\frac{1}{\sqrt2}\phi-(\sqrt2-1)-\left(1-\frac{1}{\sqrt{2}}\right)\phi^{-1}.
\end{equation}
One can solve~\eqref{FIK ode} implicitly up to an arbitrary constant $\eta$, obtaining
\begin{equation}		\label{FIK implicit}
e^{r+\eta}=\frac{\phi-1}{\big(\phi+\sqrt2-1\big)^{\sqrt2-1}}.
\end{equation}
As shown in~\cite{FIK03}, the $\mc L^2_{-1}$ soliton is complete and exists for all $r\in\mb R$. Using this fact and the positivity of $\phi$ and $\phi_r$,
 it is not difficult to see from~\eqref{FIK ode} and~\eqref{FIK implicit} that $\phi\nearrow\infty$
as $r\nearrow\infty$ and that the metric is asymptotically conical in the precise sense that
\begin{equation}
\frac{\phi_r}{\phi}\rightarrow\frac{1}{\sqrt2} \qquad\mbox{ as }\qquad \phi\nearrow\infty.
\end{equation}
Specifically, if we define $\gamma:=2^{-1/4}$, then the asymptotic cone of the blowdown soliton corresponds to
\begin{equation}	\label{FIK cone}
f=\gamma^2 s\quad\mbox{ and }\quad g=\gamma s\qquad\Rightarrow\qquad\frac{f^2}{g^2}=\gamma^2=\frac{1}{\sqrt2}.
\end{equation}
All sectional curvatures of the cone vanish except $\kappa_{23} = \frac{4(\sqrt2-1)}{s^2}$.

\section{Results of simulations}

The simulations are performed using standard numerical methods for parabolic equations: centered differences for spatial derivatives and Euler's method for time evolution.  We describe the method here in more detail.

Any function $F(t,\theta)$ is represented by the values $F^k _i$ that the function takes at points $\theta_i$ equally spaced with spacing $\Delta \theta$ and with times $t_k$ equally spaced with spacing $\Delta t$.  We use standard centered finite differences, so that $F_\theta$ and $F_{\theta \theta}$ are approximated by
\begin{equation}
{F_\theta} = {\frac {{F^k _{i+1}}-{F^k _{i-1}}} {2 \Delta \theta}}
\label{fd1}
\end{equation}
and
\begin{equation}
{F_{\theta \theta}} = {\frac {{F^k _{i+1}} + {F^k _{i-1}} - 2 {F^k _i}} {{(\Delta \theta)}^2}},
\label{fd2}
\end{equation}
respectively.

Time evolution is carried out using the Euler method, so that 
\begin{equation}
{F^{k+1} _i} = {F^k _i} + \Delta t {\partial _t} F, 
\label{Euler1}
\end{equation}
where ${\partial _t} F$ is the finite difference version of the right hand side of equation (\ref{dtA}), (\ref{dtB}) or (\ref{dtC}) with the spatial derivatives evaluated using equations (\ref{fd1}) and (\ref{fd2}).  The standard von Neumann stability analysis of equation (\ref{Euler1}) reveals that the time step must satisfy the Courant condition
\begin{equation}
\Delta t < {\textstyle {\frac 1 2}} {{({\rho _{\rm min}}\Delta \theta)}^2}.
\label{CFL1}
\end{equation}
where $\rho _{\rm min}$ is the minimum value of $\rho$.

We now show the results of a simulation with $\epsilon=0.1$.  To get an idea of how the evolution proceeds we plot $1/{\kappa_{23}}$ at $\theta=\pi/2$ as a function of time.  This curve goes to zero (i.e. ${\kappa _{23}} \to \infty$) at the final time.  Note that the curve is linear near the final time.  Thus $1/{\kappa _{23}}$ is a good proxy (up to some overall scale) for time remaining until the singularity.

\begin{figure}[H]
	\includegraphics[width=0.9\linewidth]{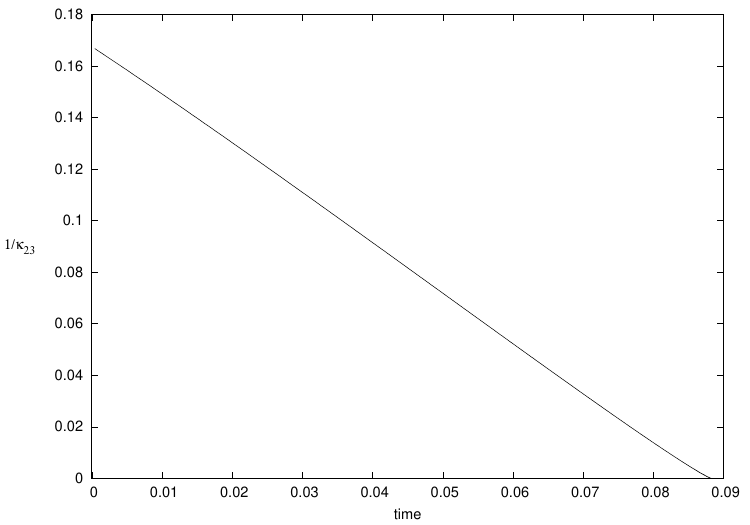}
	\caption{$1/{\kappa_{23}}$ plotted vs. time}
	\label{ftgsq}
\end{figure}

\begin{figure}[H]
	\includegraphics[width=0.9\linewidth]{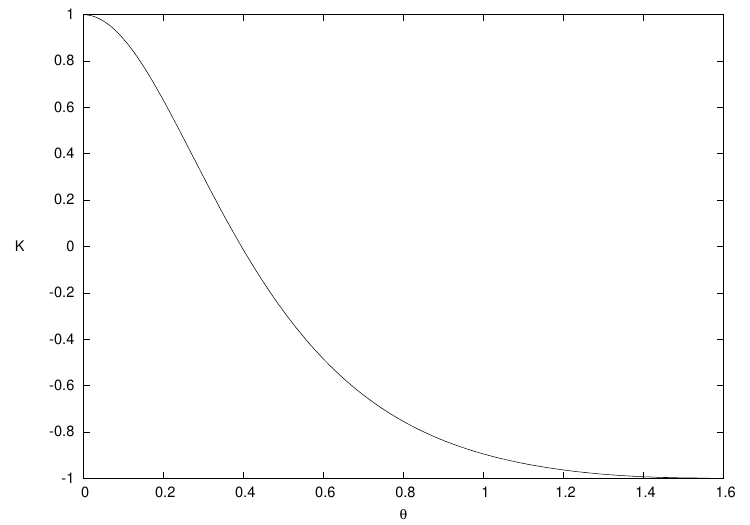}
	\caption{$K=g{g_s}/f$ plotted vs. $\theta$ near the final time}
			\label{Kahlerfig}
\end{figure}

We now want to examine the extent to which the metric becomes K\"ahler as the singularity is approached.  Recall that a K\"ahler metric has $f=g{g_s}$.  We define the quantity $K$ by $K = g{g_s}/f$.  Then a metric is K\"ahler if $K=1$, but it is also K\"ahler (but with the opposite orientation of $s$) if $K=-1$.  In Figure \ref{Kahlerfig}, we plot $K$ vs. $\theta$ at a time near the final time.  Note that near $\theta = \pi/2$ the quantity $K$ is approaching $-1$, thus indicating that the metric in this region is becoming K\"ahler.  

We would now like to know the form of the K\"ahler metric that is being approached.  Here we take a small region near $\theta=\pi/2$ and examine the behavior of the cone angle  $\gamma^2={f^2}/{g^2}$. In Figure \ref{ftgsq}, we plot $\gamma^2={f^2}/{g^2}$ vs. the length from $\pi/2$ for three times near the final time. We see that $\gamma^2$ approaches $0.707\approx 1/\sqrt{2}$ numerically, indicating local convergence of the metric cone angle $\gamma^2$ to that of the blowdown soliton; \emph{cf}. equation \eqref{FIK cone}.  What is perhaps surprising is how close one needs to get to the singularity to obtain this asymptotic behavior.  In figure (\ref{ftgsq}) the bottom curve corresponds to $1/{\kappa_{23}}=4.7 \times {{10}^{-8}}$ and has an asymptotic value for ${f^2}/{g^2}$ of 0.666. The middle curve has $1/{\kappa_{23}}=9.4 \times {{10}^{-9}}$ and has an asymptotic value for ${f^2}/{g^2}$ of 0.700.  The top curve has $1/{\kappa_{23}}=1.9 \times {{10}^{-9}}$ and has an asymptotic value for ${f^2}/{g^2}$ of 0.707.

\begin{figure}[H]
	\includegraphics[width=0.9\linewidth]{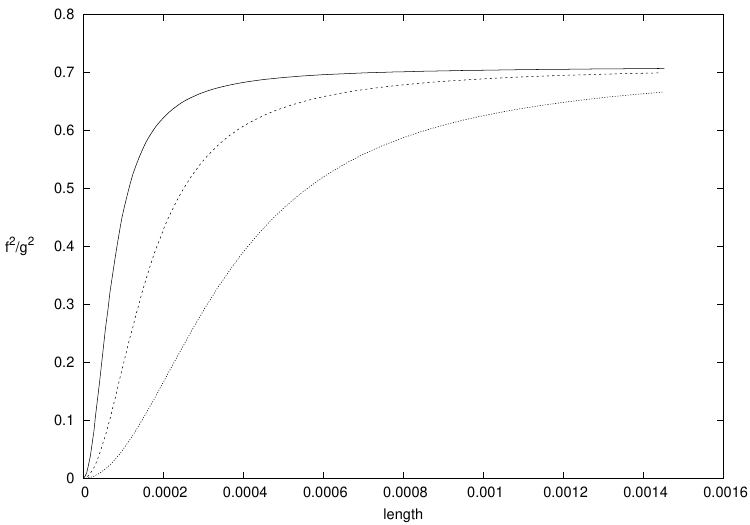}
	\caption{$\gamma^2 = f^2/g^2$ plotted vs. length near the final time}
	\label{ftgsq}
\end{figure}

\section{Directions for future work}

Our results here provide evidence in favor of the conjecture that unstable perturbations of the Fubini--Study Einstein metric --- perturbations which are conformal and not K\"ahler ---
develop finite-time local singularities modeled by the K\"ahler blowdown soliton $\mc L_{-1}^2$ discovered in~\cite{FIK03}, but with the opposite complex structure.

Short of proving the full conjecture, a useful next step would be to derive formal matched asymptotics that describe how parabolic dilations at these singularities approach
$\mc L_{-1}^2$, analogous to the matched asymptotics formally derived in~\cite{AIK11} for solutions that approach the Bryant soliton. However, because $\mc L_{-1}^2$ is
only known in the implicit form~\eqref{FIK implicit}, a somewhat more approachable next step would be to develop formal matched asymptotics for parabolic dilations that
approach its asymptotic cone $(f,g)=(\gamma^2s,\gamma s)$, where $\gamma=2^{-1/4}$, as described in~\eqref{FIK cone}.



\bibliography{fik_num}

\end{document}